**Gromov-Wasserstein and optimal transport: from assignment problems to probabilistic numeric**

Iman Seyedi[1,*], Antonio Candelieri[2], Enza Messina[1], and Francesco Archetti[1]

[1] University of Milano-Bicocca, Department of Computer Science Systems and Communication, Italy; seyediman.seyedi@unimib.it, enza.messina@unimib.it, and francesco.archetti@unimib.it
[2] University of Milano-Bicocca, Department of Economics Management and Statistics, Italy; antonio.candelieri@unimib.it
* Correspondence: seyediman.seyedi@unimib.it

**Abstract**

The assignment problem, a cornerstone of operations research, seeks an optimal one-to-one mapping between agents and tasks to minimize total cost. This work traces its evolution from classical formulations and algorithms to modern optimal transport (OT) theory, positioning the Quadratic Assignment Problem (QAP) and related structural matching tasks within this framework. We connect the linear assignment problem to Monge's transport problem, Kantorovich's relaxation, and Wasserstein distances, then extend to cases where source and target lie in different metric-measure spaces—requiring Gromov–Wasserstein (GW) distances. GW formulations, including the fused GW variant that integrates structural and feature information, naturally address QAP-like problems by optimizing alignment based on both intra-domain distances and cross-domain attributes. Applications include graph matching, keypoint correspondence, and feature-based assignments. We present exact solvers, Genetic Algorithms (GA), and multiple GW variants—including a proposed multi-initialization strategy (GW_MultiInit) that mitigates the risk of getting stuck in local optima alongside entropic Sinkhorn-based approximations and fused GW. Computational experiments on capacitated QAP instances show that GW_MultiInit consistently achieves near-optimal solutions and scales efficiently to large problems where exact methods become impractical, while parameterized EGW and FGW variants provide flexible trade-offs between accuracy and runtime. Our findings provide theoretical foundations, computational insights, and practical guidelines for applying OT and GW methods to QAP and other real-world matching problems, such as those in machine learning and logistics.

1. **Introduction**

The assignment problem is a key challenge in operations research. It is used in workforce allocation, supply chain management, machine learning, and computer vision [1], [2], [3]. The goal is to find the best one-to-one matching between two equal sets, usually agents and tasks. Each agent must be assigned to one task, and each task to one agent, with the total cost as low as possible [4]. The common way to solve this has been combinatorial optimization. The Hungarian algorithm has been the main method for years [5]. Kuhn first developed it, and Munkres refined it. It gives an exact solution in polynomial time, with O(n³) complexity for $n \times n$ assignment matrices [4], [6]. But in modern use cases, the classical approach shows limits. Many problems now involve uncertain demand, different capacity limits, and complex structural rules [7]. For example, a logistics company may need to assign vehicles to service regions. Demand forecasts might be uncertain, vehicles may not have the same capacity, and location relationships can matter as much as assignment costs. Classical assignment formulations struggle to capture such nuanced requirements [8]. Classical assignment methods work well in clear, structured cases. But they struggle when demands are uncertain, capacities are different, or when there are dependencies beyond simple pairwise costs [7], [8]. The Quadratic Assignment Problem (QAP) [9] is a generalization. It adds pairwise interaction costs between assignments. This makes it useful for problems like facility layout, graph matching, and spatial arrangements. QAP is NP-hard [10]and becomes intractable for large-scale instances, motivating the need for approximation techniques and structure-exploiting algorithms. An early result by [11] proves that even a single negative eigenvalue makes the QAP NP-hard. More recently, [12] provides further insight into the non-convex nature of the Gromov–Wasserstein (GW) distance and, in Theorem 4, establishes the NP-hardness of computing the GW distance between finite spaces for any instance of the input data. While in a few specific cases of metric-measure spaces the GW distance can be determined analytically, in the general case of finite spaces, it remains NP-hard.

The evolution from discrete assignment problems to optimal transport theory marks a major shift that helps overcome these limits [13]. Optimal transport comes from Monge's 1781 problem of moving earth with minimal effort. It gives a broad mathematical framework that extends assignment problems. It also provides tools to deal with uncertainty, structural constraints, and continuous distributions [14], [15]. This framework transforms the discrete matching problem into a measure-theoretic optimization over transport plans, enabling mass splitting, capacity constraints, and probabilistic formulations that are intractable in classical settings [16]. Contemporary applications further motivate this theoretical evolution. In machine learning, optimal transport is now a key way to compare and work with probability distributions. Surveys show its impact in areas like generative modeling, domain adaptation, and neural network optimization [17], [18]. The Wasserstein distance, which comes from optimal transport, is now a central metric in deep learning. It is widely used to train generative adversarial networks (GANs) and to measure differences between distributions [19], [20]. Recent work has taken optimal transport into new fields. It has been applied to single-cell genomics, federated learning, and large-scale neural network training [17], [21]. This shows its flexibility beyond standard machine learning tasks. In computer vision, transport-based methods work well for image registration and shape matching. Classical assignment methods cannot handle the continuous changes in images [22]. In supply chain optimization, companies face uncertain demand and multimodal transport networks, scenarios where traditional assignment methods fail but optimal transport methods thrive [23].

The computational landscape has undergone significant evolution. Exact assignment algorithms like the Hungarian method work in polynomial time for small to medium problems. But they become too costly when the scale grows to thousands of agents and tasks [24]. Modern optimal transport algorithms, especially entropy-regularized ones like the Sinkhorn algorithm, give near-optimal solutions with big computational gains [25]. The Sinkhorn algorithm uses entropic regularization of the Kantorovich formulation. It turns the linear programming problem into matrix scaling steps that are easy to parallelize and run well on GPUs [25], [26]. Scalability studies show that modern optimal transport methods can handle datasets with millions of samples. They keep efficiency by using neural optimal transport, mini-batch methods, and distributed computing frameworks [17], [18]. Tests show that Sinkhorn-based methods can be up to 100× faster while keeping solution quality close to exact methods [27], [28]. This efficiency comes from entropic regularization, which smooths the discrete optimal transport problem and makes it easier for iterative algorithms [29]. The regularization parameter gives a clear trade-off between speed and accuracy, so users can adjust it for their needs [26]. Newer methods now use neural optimal transport. These rely on deep learning to solve transport problems end-to-end, making it possible to work with high-dimensional data that were once too hard [30], [31]. The Sinkhorn algorithm's matrix-based operations leverage modern hardware architectures effectively, making it particularly suitable for large-scale applications in machine learning and data science [27].

The present paper asks a key question: How can we connect classical assignment problems with modern optimal transport theory to build a single framework that uses the efficiency of transport methods while keeping the clarity and guarantees of classical approaches? We demonstrate that quadratic assignment problems appear as special cases of optimal transport. This perspective opens new ways to handle complex real-world constraints that traditional methods cannot address [29]. By establishing this theoretical connection and providing a comprehensive computational analysis, we show how optimal transport provides both a generalization of classical assignment problems and a practical computational framework for modern large-scale applications.

Optimal Transport (OT) provides a powerful tool to compare probability measures by minimizing a transport cost between them. However, many real-world data sources—such as images, shapes, graphs, or sequences—do not consist of simple points in a shared space. Instead, they are structured objects with internal relationships, like pairwise distances between key-points in an image or connections in a graph. In these cases, standard OT is not well-suited, because it assumes that both distributions lie in the same metric space and compares points directly. Moreover, classical OT is sensitive to isometric transformations such as rotations and translations, which might be irrelevant in some problems. For example, two images may represent the same object but with different orientations.

GW distances address this problem. The GW distance, first introduced by [32], has attracted significant attention for its ability to compare structured objects represented as metric measure spaces. This makes GW useful when absolute positions are not meaningful, such as in graph matching, keypoint alignment, or cross-domain comparison. For a detailed survey on the geometry of GW distances, see [33]. However, as the scope of applications has expanded, several limitations of classical GW have motivated new formulations and improvements. To apply GW in broader contexts, several extensions have been proposed. Fused GW combines structure-based and feature-based costs [33]. Partial GW [34] allows only a fraction of the total mass to be transported, which is useful when datasets have outliers or partial correspondence. Another important contribution focuses on computational efficiency. Zhang, Wang, et al. [35] present acceleration methods for computing the gradient of the entropic GW distance. By exploiting the structure of distance matrices and applying dynamic programming, they significantly reduce the cost of matrix multiplications, a key bottleneck in GW solvers. Similarly, [36] studies the stability and algorithmic properties of the entropic GW distance. They provide a variational representation of the regularized problem and propose a convergent algorithm that is much more tractable than solving the original NP-hard GW formulation. Peyré et al. [37] proposed an entropically regularized algorithm that fixes the number of support points and the barycenter measure, making the problem tractable. This idea is used in applications like averaging graphs or shapes. Many problems also require multi-input matching. For example, clustering [38], graph partitioning [33], and multi-graph alignment—require simultaneous matching of more than two inputs. One common approach is to use a reference measure (often the GW barycenter) and compute bi-marginal transports. Alternatively, [39] introduced a multi-marginal extension of GW, generalizing multi-marginal OT to the structural setting. One major line of research concerns the unbalanced GW distance, which relaxes the constraint that both measures must be probability distributions. The work of [40] introduces two unbalanced GW formulations, allowing the comparison of metric spaces equipped with arbitrary positive measures. These formulations are based on a conic relaxation and lead to more general and robust frameworks, especially in cases involving missing data or outliers. More broadly, the review by [41] outlines how unbalanced OT, entropic regularization, and GW can be combined into scalable and flexible loss functions for modern data science tasks. Sliced GW approximates GW using random 1D projections [42]. Extensions of GW have also been developed to handle multiple distributions. In particular, Peyré et al. [43] introduce a method for Gromov-Wasserstein barycenters, which enables the averaging of distance or kernel matrices. Their method minimizes distortion through probabilistic mappings and is formulated as a regularized version of the GW problem. This is useful in tasks like clustering and structure summarization, where an average shape or graph is needed. While GW is useful for finding soft correspondences between two metric spaces (i.e., probabilistic couplings), certain applications require exact matchings or registrations. For this purpose, [44] propose the Gromov-Monge (GM) distance, which modifies GW by replacing soft correspondences with hard (Monge) mappings. Although more restrictive, GM allows for sharper structural alignment when needed. For scalability in large-scale problems, Sliced Gromov-Wasserstein (SGW) has been proposed by [42]. SGW reduces the complexity of GW by projecting the metric spaces to one-dimensional structures and applying closed-form solutions related to the quadratic assignment problem. Its computational cost of $O(n^2)$ makes it suitable for problems with large numbers of samples while maintaining a strong connection to the original GW formulation.

Together, these works demonstrate the growing versatility of GW and its variants. They provide theoretical tools and practical algorithms for handling diverse real-world tasks, including data with unbalanced masses, structural noise, and large-scale complexity. In addition to pairwise matching, there has been work on special cases. For instance, GW distances have been studied for Gaussian measures and have been extended through related concepts such as C-optimal transport [45], [46] as a generalization of BO distances and its unbalanced variants [43], [47]. Another important topic is the GW barycenter problem. The goal is to find a central representative among multiple structured inputs [39]. Recent work has looked at ways to compute distances and alignments between unaligned graphs, especially for Gromov–Wasserstein-type objectives. In [48], the authors suggest a method to find both the distance and the transportation plan between two graphs. They do this by finding the best permutation using stochastic gradient descent (SGD). Their Wasserstein structural transport (WST) cost function is differentiable, so it works with gradient-based optimization. The algorithm is non-convex, so it can get stuck in local minima. To handle this, they reformulate the objective as an expected value over a parameter space. This helps avoid bad local optima. A key part of their approach is using Bayesian exploration, which gives the algorithm global convergence properties. More generally, Bayesian optimization (BO) is a useful tool for handling non-convex permutation problems [49], [50]. In the Bayesian Optimization for Permutation Spaces (BOPS) framework, a Gaussian process surrogate model is defined over the permutation space using the Mallows kernel, which measures similarity between permutations via the exponential negative of the number of discordant pairs[51]. Building on this, [52] introduce a sorting-based kernel (Merge kernel) and its associated optimization method, MERGE-BO, for efficient BO in permutation spaces. MERGE-BO has been evaluated on both synthetic benchmarks and real-world applications, demonstrating competitive performance.

The present paper builds on these ideas by studying GW and its variants—Entropic and fused GW—for solving quadratic assignment problems between featurized sources and targets. Such problems arise in real data, where absolute distances are often unreliable but structural similarity remains important. We compare exact and regularized solvers and provide computational results on synthetic and real tasks. Our main contribution is a unified framework that rigorously connects classical assignment problems to optimal transport theory, showing they emerge as a special case of the Kantorovich formulation when transport costs correspond to discrete assignment matrices [53], [54]. We trace the evolution of computational methods from exact linear programming solutions to entropy-regularized approximations and modern neural optimal transport approaches, with complexity analysis showing the shift from O(n³) Hungarian algorithm complexity to O(n²) Sinkhorn iterations while preserving theoretical approximation guarantees [26]. Computationally, we compare exact and approximate methods across multiple problem scales, addressing discrete-discrete, semi-discrete, and continuous-continuous transport regimes. We implement advanced extensions, including Wasserstein Barycenters, Gromov-Wasserstein distances, and unbalanced transport. Using the Python Optimal Transport (POT) library [27], we provide a practical implementation framework with algorithmic descriptions, parameter tuning guidelines, and scalability analysis, demonstrating advantages in applications spanning machine learning, computer vision, and logistics.

This paper is structured as follows. Section 2 presents the transition from classical assignment problems to the optimal transport framework, establishing the theoretical connections between discrete formulations, Wasserstein distances, and the Gromov–Wasserstein setting. Section 3 presents the computational methods and algorithms. It focuses on Wasserstein distance and its variants, and explains the Gromov–Wasserstein formulation for comparing structured data. Section 4 discusses key extensions of GW, including entropic Gromov–Wasserstein and fused Gromov–Wasserstein along with the Multiple Initialization Strategy for Gromov-Wasserstein Optimization. Section 5 presents the computational results, including synthetic benchmarks and real-world assignment tasks, together with scalability and sensitivity analyses. Finally, Section 6 gives the conclusion, talks about the limitations, and outlines directions for future research.

## 2. From Classical Assignment to Optimal Transport

### 2.1 The Classical Assignment Problem

The assignment problem represents one of the fundamental and well-studied optimization problems in operations research and combinatorial optimization [55]. It addresses the challenge of finding an optimal one-to-one matching between two finite sets of equal cardinalities while minimizing total cost or maximizing total benefit. The problem's mathematical elegance, combined with its broad practical applicability, has made it a cornerstone of optimization theory and a natural stepping stone towards more complex transport formulations.

Consider two finite sets: a set of n agents $A = \{a_1, a_2, \ldots, a_n\}$ and a set of n tasks $T = \{t_1, t_2, \ldots, t_n\}$. Let $C \in \mathbb{R}^{n \times n}$ be the cost matrix where $C[i][j]$ represents the cost of assigning agent $a_i$ to task $t_j$. The classical assignment problem seeks to find a bijective mapping $\varphi: A \to T$ that minimizes the total assignment cost.

The problem can be formulated as the following integer linear program:

$$Minimize \sum_{i=1}^{n}\sum_{j=1}^{n} C[i][j] \cdot x[i][j] \tag{1}$$

Subject to:

$\sum_{j=1}^{n} x[i][j] = 1$   for all i (Each agent is assigned to exactly one task)   (2)

$\sum_{i=1}^{n} x[i][j] = 1$   for all j (Each task is assigned to exactly one agent)   (3)

$x[i][j] \in \{0,1\}$   for all i, j   (4)

Here, x[i][j] is a binary decision variable that equals 1 if agent i is assigned to task j, and 0 otherwise. The constraints ensure that each agent is assigned to exactly one task, and each task is assigned to exactly one agent, creating a perfect match.

### 2.2 Monge's Transport Problem

The historical foundations of optimal transport trace back to Gaspard Monge's seminal 1781 problem, which addressed the practical question of moving piles of material (e.g., sand) to fill holes with minimal effort. Monge's formulation laid the groundwork for modern assignment and optimal transport theory, providing a mathematical framework for mapping one distribution of mass to another with minimal cost [13]. Let $\mu = \sum_{i=1}^{n} a_i \delta_{x_i}$ and $\nu = \sum_{j=1}^{n} b_j \delta_{y_j}$ be probability measures supported on metric spaces $(\mathcal{X}, d_\mathcal{X})$ and $(\mathcal{Y}, d_\mathcal{Y})$, respectively. The discrete Monge problem seeks a measurable map $T: X \to Y$ that pushes $\mu$ forward to $\nu$ (denoted $T \# \mu = \nu$ while minimizing the total transport cost:

$$\min_{T} \sum_{i=1}^{n} c(x_i, T(x_i)) \qquad \text{subject to} \quad \sum_{i} a_i = \sum_{j} b_j, \quad \forall j \in \{1, \dots, m\}. \tag{5}$$

In the special case n=m and uniform weights $a_i = b_j = \frac{1}{n}$, T becomes a bijection, reducing the Monge problem to the classical assignment problem, where the cost matrix $C_{i,j} = c(x_i, y_j)$

**Continuous Monge Problem:**

$$\min_T \int_X c(x, T(x)) \, d\mu(x) \quad \text{subject to } T\#\mu = \nu \tag{6}$$

where $c(x,y)$ represents the cost of transporting a unit mass from $x$ to $y$. For measures with densities $f(x)$ and $g(y)$, this becomes:

$$\min_T \int_X c(x, T(x)) f(x) \, dx \quad \text{subject to } T\#\mu = \nu \tag{7}$$

Discrete Formulation: For discrete distributions with finite support, let:

- Source points: $X = \{x_1, x_2, \ldots, x_n\}$ with masses $a = (a_1, a_2, \ldots, a_m)$
- Target points: $Y = \{y_1, y_2, \ldots, y_n\}$ with masses $b = (b_1, b_2, \ldots, b_n)$
- Cost function: $c(x_i, y_j)$ representing the transport cost from $x_i$ to $y_j$

The classical assignment problem emerges as a special case of the discrete Monge transport problem under specific conditions:

Table 1. Key Conceptual Differences between Monge and Assignment Problem

| Aspect | Monge Transport | Assignment Problem |
|---|---|---|
| Source points | Arbitrary locations with masses | Agents (unit masses) |
| Target points | Arbitrary locations with demands | Tasks (unit demands) |
| Cost function | General transport cost c(x,y) | Assignment cost matrix C[i][j] |
| Mass distribution | General positive masses | Uniform: $a_i = b_j = 1$ |
| Objective | Minimize total transport cost | Minimize total assignment cost |
| Constraint | No mass splitting | One-to-one assignment |
| Solution structure | General transport map | Permutation matrix |

When source and target distributions are uniform ($\sum_i a_i = \sum_j b_j = 1$ for all i, j) and the number of agents equals the number of tasks (m = n), the Monge problem reduces to:

Minimize $\sum_{i=1}^{n} C[i][\sigma(i)]$ (8)

where $\sigma: \{1, \ldots, n\} \rightarrow \{1, \ldots, n\}$ is a permutation representing the assignment strategy.

### 2.3 Kantorovich Relaxation

Leonid Kantorovich's 1942 relaxation of Monge's problem was a key step in optimal transport theory, with wide theoretical and computational use [13]. Unlike Monge's deterministic transport maps $T: X \rightarrow Y$, Kantorovich's approach uses a transport plan. This allows mass from one source to be split and sent to multiple targets. It creates a probabilistic link between source and target measures, expands the set of possible solutions, and ensures solutions exist under mild conditions. The Kantorovich formulation replaces Monge's deterministic transport maps with probabilistic transport plans. For probability measures $\mu$ and $\nu$, a transport plan $\pi$ is a probability measure on the product space $\mathcal{X} \times \mathcal{Y}$ with marginals $\mu$ and $\nu$.

**Continuous Kantorovich Problem:** Let $\Pi(\mu, \nu)$ denote the set of all probability measures on $\mathcal{X} \times \mathcal{Y}$ with marginals $\mu$ and $\nu$. The Kantorovich optimal transport problem is:

$$\min_{\pi \in \Pi(\mu,\nu)} \int_{\{\mathcal{X} \times \mathcal{Y}\}} c(x, y) \, d\pi(x, y) \qquad (9)$$

where $c(x,y)$ denotes the cost of moving a unit mass from $x$ to $y$. This formulation reduces to Monge's problem when $\pi$ is concentrated on the graph of a measurable map $T$.

**Discrete Kantorovich Problem**: For discrete distributions with source masses $a = (a_1, \ldots, a_m)$ and target masses $b = (b_1, \ldots, b_n)$, $a_i, b_j \geq 0$, $\sum_i a_i = \sum_j b_j = 1$ let the transport plan become a matrix $T \in \mathbb{R}_+^{\{m \times n\}}$ where $T[i][j]$ represents the mass transported from source i to target j.

Minimize $\sum_{i=1}^{m} \sum_{j=1}^{n} C[i][j] \cdot T[i][j]$ (10)

Subject to:

$\sum_{j=1}^{n} T[i][j] = a_i$  for all $i = 1, \ldots, m$ (source constraints) (11)

$\sum_{i=1}^{m} T[i][j] = b_j$  for all $j = 1, \ldots, n$ (target constraints) (12)

$T[i][j] \geq 0$  for all i, j  (non-negativity) (13)

The assignment problem emerges as a special case of the Kantorovich formulation where the optimal transport plan is restricted to permutation matrices:

Table 2. Key Conceptual Differences between Kantorovich and Assignment Problem

| Feature | Kantorovich Transport | Assignment Problem |
| --- | --- | --- |
| **Mass splitting** | Allowed | Not allowed |
| **Transport plan** | General matrix $T[i][j] \in \mathbb{R}_+$ | Permutation matrix $T[i][j] \in \{0,1\}$ |
| **Feasible region** | Convex polytope | Vertices of assignment polytope |
| **Capacity constraints** | Naturally incorporated | Fixed unit capacities |
| **Optimization type** | Convex (linear programming) | Combinatorial |
| **Solution uniqueness** | Generally not unique | Unique under non-degeneracy |

When source and target distributions are uniform ($a_i = b_j = 1$), the assignment problem corresponds to finding extreme points of the Kantorovich polytope, which are precisely the permutation matrices.

### 2.3.1 Fundamental Advantages Over Monge Formulation

**Mass Splitting Permitted**: The Kantorovich relaxation allows mass at each source to be distributed among multiple targets, dramatically expanding the feasible solution space and ensuring existence of optimal solutions under mild conditions.

**Convex Optimization Structure**: The discrete Kantorovich problem is a linear program, making it amenable to efficient solution methods including simplex algorithms, interior point methods, and specialized network flow algorithms.

**Guaranteed Existence**: Under standard assumptions (compact support, continuous cost functions), optimal transport plans always exist, eliminating the existence issues that plague Monge problems.

**Dual Formulation**: The Kantorovich formulation admits a natural dual problem that provides economic interpretation and computational advantages through complementary slackness conditions.

### 2.4. Quadratic Assignment Problem

The QAP is a classic combinatorial optimization problem. It generalizes the linear assignment problem by adding pairwise interaction costs between assignments [9], [10]. Koopmans and Beckmann first formulated it for facility layout [9]. The QAP models situations where the cost of assigning an agent to a task depends not only on that assignment but also on the assignments of other agents. This quadratic dependency makes QAP a natural framework for facility location, electronic circuit design, graph matching, and structural pattern recognition [56].

Let $F \in \mathbb{R}^{n \times n}$ be the *flow matrix*, where $F[i][k]$ represents the amount of interaction (flow) between facility i and facility k. $D \in \mathbb{R}^{n \times n}$ be the *distance matrix*, where $D[j][l]$ represents the distance (or dissimilarity) between location j and location l.

The QAP seeks a bijective assignment $\sigma: \{1, \dots, n\} \to \{1, \dots, n\}$ that minimizes the total interaction cost:

$$\min_{\sigma \in S_n} \sum_{i=1}^{n} \sum_{k=1}^{n} F[i][k] \cdot D[\sigma(i)][\sigma(k)] + \sum_{i=1}^{n} C[i][\sigma(i)] \tag{12}$$

where:

- The first term captures *quadratic costs* due to interactions between pairs of facilities and the corresponding distances between their assigned locations.
- The optional second term $C[i][\sigma(i)]$ represents *linear assignment costs* that may be present in hybrid formulations.

The binary integer programming form of QAP uses decision variables $x_{ij} \in \{0,1\}$ indicating whether facility i is assigned to location j:

$$\min_x \sum_{i=1}^{n} \sum_{k=1}^{n} \sum_{j=1}^{n} \sum_{l=1}^{n} F[i][k] \cdot D[j][l] \cdot x_{ij} \cdot x_{kl} + \sum_{i=1}^{n} \sum_{j=1}^{n} C_{ij} \cdot x_{ij} \tag{13}$$

Subject to:

$\sum_{j=1}^{n} x_{ij} = 1 \qquad \forall i \text{(each facility is assigned to one location)}$ (14)

$\sum_{i=1}^{n} x_{ij} = 1 \qquad \forall j \text{(each location hosts one facility)}$ (15)

$x_{ij} \in \{0,1\} \qquad \forall i,j$ (16)

This quadratic objective makes the problem NP-hard [57], rendering it intractable for $n \gtrsim 30$ in most cases.

### 2.4.1 Relationship between optimal transport and QAP

The QAP can be viewed as a structural matching problem, where the goal is to minimize discrepancies between two pairwise relationship matrices under a permutation. This aligns naturally with the GW optimal transport framework [32], [33], which compares metric-measure spaces by aligning their intra-domain distance matrices. Specifically, the quadratic term

$$\sum_{i,k,j,l} F[i][k] \cdot D[j][l] \cdot x_{ij} \cdot x_{kl} \tag{17}$$

is equivalent in form to the GW cost function when F and D are interpreted as intra-domain distances (or similarities) in the source and target spaces.

### 2.4.2 Capacitated Quadratic Assignment Problem

The Capacitated Quadratic Assignment Problem (CQAP) extends the classical QAP by introducing capacity constraints on facilities (agents) and demand requirements on locations (tasks). Unlike the standard QAP, where each facility is assigned to exactly one location, the CQAP allows facilities to serve multiple locations, subject to capacity limits, while locations must satisfy their demand requirements. This formulation is especially relevant in logistics, supply chain design, and telecommunication networks [56], [58].

Let:

- $u_i$: capacity of facility i.
- $d_j$: demand at location j.

The CQAP objective is:

$$\min_{x} \sum_{i=1}^{n}\sum_{k=1}^{n}\sum_{j=1}^{m}\sum_{l=1}^{m} F[i][k] \cdot D[j][l] \cdot x_{ij} \cdot x_{kl} + \sum_{i=1}^{n}\sum_{j=1}^{m} C_{ij} \cdot x_{ij} \tag{13}$$

subject to capacity and demand constraints:

$$\sum_{j=1}^{m} d_j x_{ij} \leq u_i, \forall i = 1, \dots, n \quad \text{(Each facility cannot exceed its capacity)} \tag{14}$$

$$\sum_{i=1}^{n} u_i x_{ij} \geq d_j, \forall j = 1, \dots, m \quad \text{(Each location's demand must be satisfied)} \tag{15}$$

$$x_{ij} \in \{0,1\}, \forall i,j \tag{16}$$

This capacity constraint extension makes the problem even harder computationally than standard QAP, belonging to the class of NP-hard problems. Nevertheless, it provides a realistic modeling tool for problems involving both interaction costs and heterogeneous resource constraints.

## 3. Computational Methods and Algorithms

### 3.1 Wasserstein Distances as Transport Costs

The Wasserstein distance is a key link between optimal transport theory and modern statistics. It provides a solid way to compare probability distributions and has impacted machine learning, statistics, and computational optimization [59], [60]. Unlike traditional measures like Kullback-Leibler divergence or total variation distance, Wasserstein distances follow metric rules and give meaningful geometric measures of similarity that respect the structure of the space.

### 3.1.1 Mathematical Definition and Properties

The p-Wasserstein distance between probability measures α and β supported on a metric space (Ω, d) is defined through the Kantorovich formulation [59]:

$$W_p^p(\alpha, \beta) = \inf_{\pi \in \Pi(\alpha, \beta)} \int_{\Omega \times \Omega} d^p(\omega, \omega') d\pi(\omega, \omega') \quad (17)$$

where $\Pi(\alpha, \beta)$ denotes the set of all probability measures on $\Omega \times \Omega$ with marginals α and β, and p ≥ 1. This formulation directly connects Wasserstein distances to optimal transport by interpreting the distance as the minimum expected cost of transporting mass from distribution α to distribution β under cost function $c(\omega, \omega') = d^p(\omega, \omega')$.

### 3.2.2 The 2-Wasserstein Distance and Monge Formulation

The 2-Wasserstein distance, corresponding to p = 2 with Euclidean ground distance, admits an equivalent Monge formulation that provides additional geometric insight:

$$W_2^2(\alpha, \beta) = \min_{T_{\#}\alpha = \beta} \int_{\Omega} \| T(\omega) - \omega \|_2^2 \, d\alpha(\omega) \quad (18)$$

which is intended as an optimal transport problem [29], [61]. where $T_{\#}\alpha = \beta$ denotes the push-forward condition ensuring that the transport map T transforms distribution α into distribution β. This formulation interprets the 2-Wasserstein distance as the minimum kinetic energy required to morph one distribution into another through a deterministic transport map.

The 2-Wasserstein distance can be understood as the geodesic distance in the space of probability measures equipped with the Riemannian metric induced by optimal transport. This geometric perspective has profound implications for:

4. Gradient flows in probability space
5. Wasserstein barycenters as geometric means
6. Interpolation between distributions along geodesics

2.4.3 Analytical Solutions for Gaussian Distributions

For Gaussian distributions $\alpha = \mathcal{N}(m_\alpha, \Sigma_\alpha)$ and $\beta = \mathcal{N}(m_\beta, \Sigma_\beta)$, the 2-Wasserstein distance admits a closed-form solution:

$$W_2^2(\alpha, \beta) = \| m_\alpha - m_\beta \|^2 + \mathcal{B}(\Sigma_\alpha, \Sigma_\beta)^2 \tag{19}$$

where $\mathcal{B}(\Sigma_\alpha, \Sigma_\beta)$ represents the Bures metric between positive definite covariance matrices [62], [63]:

$$\mathcal{B}(\Sigma_\alpha, \Sigma_\beta)^2 = trace\left(\Sigma_\alpha + \Sigma_\beta - 2\left(\Sigma_\alpha^{\frac{1}{2}} \Sigma_\beta \Sigma_\alpha^{\frac{1}{2}}\right)^{\frac{1}{2}}\right) \tag{20}$$

Thus, in the case of centered GDs (i.e., $m_\alpha = m_\beta = 0$), the 2-Wasserstein distance resembles the Bures metric. Moreover, if $\Sigma_\alpha$ and $\Sigma_\beta$ are diagonal, the Bures metric is the Hellinger distance while in the commutative case, that is $\Sigma_\alpha \Sigma_\beta = \Sigma_\beta \Sigma_\alpha$, the Bures metric is equal to the Frobenius norm $\| \Sigma_\alpha^{\frac{1}{2}} - \Sigma_\beta^{\frac{1}{2}} \|_{Frobenius}^2$.

When dealing with general probability distributions beyond the Gaussian case, computing the Wasserstein distance becomes much more difficult and requires solving the Monge or Kantorovich optimal transport problems. The Kantorovich formulation turns this into a linear program with $O(n^3 \log n)$ complexity for discrete measures on n points, which makes exact computation impractical for large-scale applications. This computational burden has motivated the development of entropic regularization techniques that approximate the optimal transport solution while maintaining computational tractability.

The computational tractability of Wasserstein-based methods has been revolutionized through entropic regularization approaches, most notably the introduction of Sinkhorn distances [25], [64], which enables rapid computation of optimal transport through entropic penalties. Recent complexity analyses have improved our understanding of entropic regularized algorithms for optimal transport between discrete probability measures, establishing refined convergence bounds and computational guarantees [65]. The extension to Wasserstein barycenter computation through fast Sinkhorn-based algorithms [59], [66] has established the algorithmic foundations for practical implementation while maintaining essential geometric properties with significant computational speedups.

## 3.2. Gromov Wasserstein

Optimal Transport enables computing assignments (or permutations) that account for metric similarities between samples. However, such a metric may not exist when the measures are defined on two different spaces. GW aims to define a meaningful distance between two metric spaces $G_1$ and $G_2$, each represented as a probability measure over their elements and attributes. Specifically, let

- $G_1 \sim \mu = \sum_{i=1}^{n} h_i \delta(x_i, a_i)$,
- $G_2 \sim \nu = \sum_{j=1}^{m} g_j \delta(y_j, b_j)$,

where:

- $h \in \Delta^n, g \in \Delta^m$ are discrete probability distributions (histograms)
- $x_i, y_j$ are structural embeddings (e.g., node indices),
- $a_i, b_j$ are node features,
- $\delta$ denotes the Dirac measure.

Let $\Pi(h, g)$ be the set of admissible couplings:

$$\Pi(h, g) = \{\pi \in \mathbb{R}_+^{n \times m} \mid \sum_{j=1}^{m} \pi_{ij} = h_i, \sum_{i=1}^{n} \pi_{ij} = g_j\} \qquad (21)$$

where $\pi_{ij}$ represents the mass transported from node i in $G_1$ to node j in $G_2$.

The GW distance is a strong generalization of the classic Wasserstein distance. It solves a key problem: comparing distributions that are not in the same space. The standard Wasserstein distance $W_2$ needs both distributions in a common space to compute direct pairwise distances. Gromov-Wasserstein, however, compares distributions using their internal structural relationships. It extends optimal transport to cases where distributions cannot be compared point by point [32]. Instead of computing distances between individual samples from two distributions, GW realigns the metric spaces by comparing their internal distance structures through a transport between distance matrices. The GW distance is designed to compare structured mm spaces sounds like a limitation, where the intrinsic structure of the objects is more important than their embedding in a common space.

Let:

- $C_1 \in \mathbb{R}^{n \times n}, C_2 \in \mathbb{R}^{m \times m}$ denote the pairwise structural dissimilarities (e.g., shortest-path distances) within each graph.

The Gromov-Wasserstein discrepancy is defined as:

$$GW_q(C_1, C_2; h, g) = \min_{\pi \in \Pi(h,g)} \sum_{i,j,k,l} |C_1(i,j) - C_2(k,l)|^q * \pi_{i,k} * \pi_{j,l} \tag{22}$$

This quantity compares the internal metric structure of $G_1$ and $G_2$, by measuring how well the pairwise distances in $G_1$ match those in $G_2$ under the coupling $\pi$. Unlike classical optimal transport, GW distance does not require the nodes to live in a common feature space Figure (1).

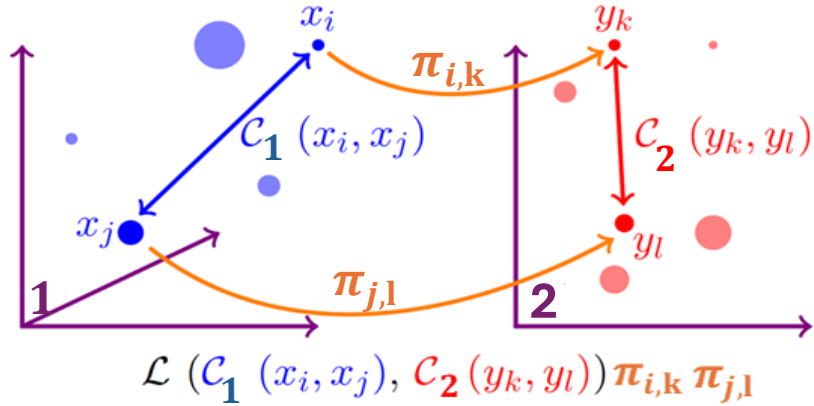

Figure 1: Illustration of GW related to Eq. (22) [67].

The assignment between two graphs (graph isomorphism) is known to be NP-hard [57]. A problem is NP-hard if the existence of a polynomial-time solution algorithm would imply $P = NP$. In molecular biology, large collaborative projects often involve gene expression measurements conducted using different machines and standards. To properly address such cases, one must consider metric measure spaces (mm-spaces), which jointly account for both the metric structure and the measure. One must consider metric measure spaces, mm-spaces, which take into account jointly the measure and the metric. Mm is a triple $\mathcal{X} = (X, d_x, \alpha)$ where $X$ is a space endowed with a metric $d_x$ and $\alpha \in M^+(X)$. A large variety of spaces are m-m: an example is a point cloud in $X = \mathbb{R}^d$ n-dimensional endowed with euclidean distance $d_x = \|\cdot - \cdot\|_2$ and the uniform distribution over the samples $\alpha = \frac{1}{N}\sum_i \delta_{x_i}$. A second is a graph where $X$ is the set of nodes, $d_x$ is the shortest path distance and $\alpha$ is the uniform probability over the graph [41].

This optimization problem is a non-convex quadratic program, making it computationally challenging. The similarity with OT stems from the use of the constraint set $\mathcal{U}(\alpha, \beta)$. the key difference is the quadratic dependence in the plan $\pi \otimes \pi$ which characterizes GW as a QAP and makes it an NP-hard problem [12], [68]. Since GW is a non-convex problem, the computations only guarantee an output $\pi$ that is the stationary point of the functional. One approach is to linearize the quadratic program [32] obtaining a sequence where one alternates between sequential updates. There exist variants of GW that are motivated by additional information. An alternative solution proposes computationally efficient lower bounds, global histograms of distances represented in [69] and tractable using OT and used as initialization.

Table 3. Key Conceptual Differences between Wasserstein Distance and Gromov-Wasserstein Distance

| Aspect | Wasserstein Distance ($W_2$) | Gromov-Wasserstein Distance |
| --- | --- | --- |
| Comparison scope | Points in the same metric space | Points in different/unrelated spaces |
| Cost computation | Direct Euclidean distance between $x_i$ and $y_j$ | Difference between internal pairwise distance structures |
| Positional requirements | Requires absolute positions | Only requires relative distances (structure-preserving) |
| Primary application | Matching with known shared geometry | Matching based on structural similarity |
| Geometric invariance | Translation and rotation sensitive | Invariant to isometric transformations |

## 4. Extensions

While the standard GW distance is effective for comparing structured data, many real-world problems require additional flexibility. To address this, several extensions of GW have been proposed. These include the Fused Gromov-Wasserstein distance, which incorporates both structural and feature information and Entropic Gromov Wasserstein. This section introduces each of these extensions and highlights their motivations and use cases.

### 4.1. Entropic Gromov Wasserstein

The GW distance involves solving a non-convex quadratic assignment problem, known to be NP-hard [35]. Traditional exact optimization techniques and convex relaxations (e.g., semidefinite programming or eigenvalue relaxations) often suffer from scalability issues, as they may require up to $O(n^4)$ variables (e.g., $n^2 \times n^2$ for convex formulations) [43], [70]. To address these computational limitations, Entropic Gromov-Wasserstein (EGW) was introduced as a scalable, differentiable approximation to the original GW problem by adding an entropy regularization term to the optimization objective [43]. Unlike classical methods that rely on rigid constraints and non-convex optimization schemes, EGW maintains the coupling constraints while providing global convergence guarantees and improved numerical stability. One of the key strengths of the entropic formulation is its flexibility: it naturally extends to structured transport problems, enabling efficient solutions for FGW and UGW. The algorithm operates in an iterative framework, where each iteration consists of two main steps [33], [35]:

- Gradient computation of the GW objective via tensor contractions, which has a computational cost of $O(n^3)$,
- Sinkhorn updates to solve the resulting entropic OT subproblem have a per-iteration cost of $O(n^2)$, but since these updates are performed repeatedly within the GW optimization loop, the overall computational complexity of the full GW algorithm remains $O(n^3)$.

This separation of the GW structure-matching step and the entropy-regularized transport update makes EGW particularly effective for medium-scale graphs. However, the cubic complexity in the number of nodes remains a bottleneck in large-scale settings. Table 4 illustrates the computational complexity and performance trade-offs of different acceleration approaches.

Table 4. The comparison of different methods for the computation of GW metric and its variants [35].

| Method | Complexity | Exact and Full-sized Plan |
|---|---|---|
| Entropic GW and its approximations | | |
| Entropic GW [43] | $O(N^3)$ | ✓ |
| S-GWL [39] | $O(N^2 \log N)$ | not exact |
| SaGroW [67] | $O(N^2(s + \log N))$ | not full-sized |
| Spar-GW [71] | $O(N^2 + s^2)$ | not full-sized |
| LR-GW [72] | $O(Nr^2d^2)$ | not exact |
| AE [73] | $O(N^2 \log N)$ | not exact |
| GW on special structures | | |
| Sliced GW [42] | $O(N^2)$ | 1D space only |
| FlowAlign [74] | $O(N^2)$ | tree only |
| FGC-GW [35] | $O(N^2)$ | ✓ |

The entropic regularized version of the optimal transport problem adds a Kullback–Leibler (KL) divergence term. The entropic GW formulation is:

$$OT_\varepsilon(h, g) = \min_{\pi \in \Pi(h,g)} \ell(\pi) + \varepsilon\, KL(\pi \parallel h \otimes g) \ , \tag{23}$$

Here:

- $\ell(\pi)$ is the standard GW loss (see Section 3.1),
- $\varepsilon > 0$ is the entropic regularization parameter,
- KL is the Kullback-Leibler divergence,
- $\Pi(h, g) = \{\pi \in \mathbb{R}_+^{n \times m} \mid \pi 1_m = h,\ \pi^T 1_n = g\}$ is the set of admissible couplings.

This regularized problem can be solved efficiently using the Sinkhorn algorithm, which performs iterative matrix scaling. Its advantages include:

- Linear convergence rate,
- GPU-accelerated matrix operations,
- Memory and time efficiency for large graphs.

Sinkhorn-type algorithms have been adapted to GW and FGW problems by relaxing the cost terms and performing iterative updates over the coupling $\pi$, often using Sinkhorn steps as subroutines (e.g., [36], [43]).

Entropic variants are especially valuable in:

- Large-scale settings, where exact solvers are infeasible,
- Soft-matching scenarios where approximate transport is acceptable,
- Unbalanced OT, where relaxed mass constraints integrate naturally with KL divergence terms.

### 4.2. Fused Gromov Wasserstein

While GW captures structural similarity, many practical graphs also come with node features (e.g., text, categories, vectors). The FGW distance extends GW by combining both feature-level and structure-level comparisons in a unified framework.

Let:

- $M_{ij} = d(a_i, b_k)$ be a feature dissimilarity cost (e.g., Euclidean or cosine distance),

- $C_1, C_2$ be the structure matrices as in Section 3.1.

We define a 4D structural distortion tensor:

$$L_{ijkl}(C_1, C_2) = | C_1(i,j) - C_2(k,l) |^q, \tag{24}$$

and formulate the FGW distance as a trade-off between node feature similarity and graph structure similarity [33]:

$$FGW_q^\alpha(\mu, \nu) = \min_{\pi \in \Pi(h,g)} \varepsilon_q^\alpha(\pi; M, C_1, C_2) \tag{25}$$

with the loss function:

$$\varepsilon_q^\alpha(\pi; M, C_1, C_2) = (1-\alpha) \sum_{i,k} \pi_{ik} d(a_i, b_k)^q + \alpha \sum_{i,j,k,l} |C_1(i,j) - C_2(k,l)|^q * \pi_{i,k} * \pi_{j,l} \tag{26}$$

where $\alpha \in [0,1]$ is a trade-off parameter:

- $\alpha = 0$: pure feature-based comparison (Wasserstein),
- $\alpha = 1$: pure structure-based comparison (Gromov-Wasserstein),
- $\alpha \in (0,1)$: fused distance accounting for both aspects.

### 4.3. Multiple Initialization Strategy for Gromov-Wasserstein Optimization

The Gromov-Wasserstein problem presents significant computational challenges because its objective function is non-convex. This creates many local optima that can trap optimization algorithms. Standard GW solvers usually start with a single initialization, often using the product of marginal distributions. But this often leads to suboptimal solutions, especially with complex cost matrix structures. To improve the exploration of the optimization space and mitigate the risk of getting stuck in local minima, we adapted an approach borrowed from non-convex optimization, called multi-initialization (GW_MultiInit), which runs the standard GW algorithm from multiple randomly generated starting points and selects the solution with the lowest objective value [40]. The specifics of the algorithm developed are provided in the following pseudo-code Figure (2):

**Input:** Cost matrices $C_1 \in \mathbb{R}^{n \times n}$, $C_2 \in \mathbb{R}^{m \times m}$, marginals $\mu \in \mathbb{R}^n$, $\nu \in \mathbb{R}^m$, Multi start ($T$)

**Output:** Optimal coupling $P^* \in \mathbb{R}^{n \times m}$

**Initialize:** $f^* \leftarrow +\infty$, $P^* \leftarrow \emptyset$

1. Try default GW: $P_0 \leftarrow GW(C_1, C_2, \mu, \nu)$, $f_0 \leftarrow \mathcal{L}_{GW}(P_0)$ **if** $f_0 < f^*$ then $f^* \leftarrow f_0, P^* \leftarrow P_0$

2. for $t = 1, \ldots, T$ do

   - Generate random coupling: $G_t \leftarrow Uniform(0,1)^{n \times m} + \epsilon$
   - repeat:
     - $G_t \leftarrow diag(\mu \oslash G_t 1_m) G_t$
     - $G_t \leftarrow G_t diag(\nu \oslash G_t^T 1_n)$ until $\| G_t 1_m - \mu \|_\infty < \delta$ and $\| G_t^T 1_n - \nu \|_\infty < \delta$
   - Solve GW with random init: $P_t \leftarrow GW(C_1, C_2, \mu, \nu, G_t)$
   - Evaluate: $f_t \leftarrow \mathcal{L}_{GW}(P_t)$ if $f_t < f^*$ then $f^* \leftarrow f_t, P^* \leftarrow P_t$

3. return $P^*$

Where: $\mathcal{L}_{GW}(P) = \sum_{i,j,k,\ell} (C_1^{ij} - C_2^{k\ell})^2 P_{ik} P_{j\ell}$, $\delta = 10^{-12}$

Figure 2: Pseudo-code for Gromov-Wasserstein with Multiple Random Initializations (GW-MultiInit)

The algorithm begins by solving the standard GW problem with default initialization, then performs T additional trials where each trial starts from a randomly generated valid coupling matrix. For each random initialization, we first generate a uniform random matrix and project it onto the coupling polytope $\Pi(\mu, \nu)$ using Sinkhorn iterations to ensure marginal constraints are satisfied. The GW algorithm is then executed from this random starting point, and the solution with the lowest objective value $\mathcal{L}_{GW}(P) = \sum_{i,j,k,\ell} (C_1^{ij} - C_2^{k\ell})^2 P_{ik} P_{j\ell}$ across all trials is returned as the final result. This multi-initialization strategy leverages the fact that different starting points can lead the iterative GW solver to different local optima, thereby increasing the likelihood of finding a globally optimal or near-optimal solution to the challenging non-convex Gromov-Wasserstein matching problem.

### 4.4. Genetic Algorithm

Genetic algorithms (GAs) have been widely applied to the QAP, a classic NP-hard combinatorial optimization problem where the objective is to assign facilities to locations to minimize a quadratic cost function based on flows and distances [75], [76]. In a GA framework, candidate assignments are encoded as chromosomes, often using permutation-based representations to respect the one-to-one assignment constraint. At the same time, crossover and mutation operators are adapted to maintain feasibility and effectively explore the search space. Selection mechanisms such as tournament or roulette-wheel selection guide the evolution toward high-quality solutions, and problem-specific operators, such as swap or inversion mutations, are frequently employed to accelerate convergence [77], [78]. Although GAs do not guarantee optimality, they are well-suited for large instances where exact methods become computationally infeasible, and hybrid approaches that combine GAs with local search or problem-specific heuristics have shown improved performance in both solution quality and runtime [79].

5. **Computational results:**

In our computational experiments, we address the CQAP through the lens of GW optimization, a framework originally developed for comparing structured distributions in different metric spaces. The classical CQAP seeks to assign tasks to agents in a way that minimizes a quadratic cost function, allowing where required capacity and demand constraints. The GW distance naturally aligns with this objective, as it models the assignment problem by comparing pairwise relational structures (distances) between agents and tasks. By incorporating capacities and demands into the GW formulation, we capture the essence of CQAP, where not only the matching quality (in terms of relational preservation) but also the feasibility of assignments (in terms of marginal constraints) is critical. We conduct a systematic evaluation across problem sizes—small, medium, and large—using several methods, including exact quadratic solvers, standard GW, GW_MultiInit, entropic GW (EGW), and Fused GW (FGW), each providing a different balance between computational efficiency and objective quality. Our goal is to assess how well these methods approximate the CQAP solution under scalability constraints, and to identify which variants are most suitable for practical deployment depending on the problem size and resource availability.

Standard GW solvers usually start with a single initialization, often using the product of marginal distributions. But this often leads to suboptimal solutions, especially with complex cost matrix structures. To fix this, we propose a multiple initialization strategy (GW_MultiInit). It runs the standard GW algorithm from several randomly generated starting points and picks the solution with the lowest objective value. To evaluate different methods on the CQAP, we create a set of synthetic test instances divided by problem size: Small, Medium, and Large (Table 5). Each case specifies the number of agents and tasks. Agent and task positions are randomly generated in 2D space using uniform sampling from the range $[0,10]^2$. Capacities and demands are then assigned as random integers:

Table 5. Test Problem Instances

| Size | Test ID | # Agents | # Tasks | Total Mass |
|---|---|---|---|---|
| Small | S1 | 3 | 3 | 9 |
| | S2 | 4 | 4 | 10 |
| | S3 | 5 | 6 | 17 |
| | S4 | 6 | 5 | 21 |

| Size | Test ID | # Agents | # Tasks | Total Mass |
|---|---|---|---|---|
| Medium | M1 | 10 | 10 | 31 |
|  | M2 | 12 | 14 | 36 |
|  | M3 | 15 | 12 | 43 |
|  | M4 | 20 | 20 | 62 |
| Large | L1 | 30 | 30 | 86 |
|  | L2 | 40 | 50 | 118 |
|  | L3 | 50 | 40 | 146 |
|  | L4 | 60 | 60 | 176 |
|  | L5 | 100 | 100 | 294 |

Table 6 reports the results of applying various methods to solve a range of CQAP instances of increasing size. Each instance is characterized by the number of agents and tasks, along with the corresponding total mass (sum of capacities and demands). For each case, the method(s) achieving the lowest objective value are listed as the *Best Method*, along with the achieved *Objective* value and corresponding *Runtime*. The exact method works best on small instances but becomes very slow for larger ones. Gromov-Wasserstein methods—especially the multi-initialization variant (GW_MultiInit)—consistently give high-quality approximations with much shorter runtimes. This shows they scale well and are practical for large CQAP problems. All experiments were run on a laptop with a 13th Gen Intel Core i7-1355U processor (1.70 GHz) and 16 GB RAM, using a 64-bit Windows operating system. The full implementation of all algorithms and the experimental setup is available on GitHub at https://github.com/iman-ie/GW_CQAP

for reproducibility and further research.Table 6: Performance Summary of CQAP Instances

| Instance (Agents×Tasks) | Best Method | Objective |
|---|---|---|
| 3×3 | Exact, GW_MultiInit | 981.6324 |
| 4×4 | Exact, GW_MultiInit, EGW | 821.4323 |
| 5×6 | Exact | 1420.1807 |
| 6×5 | Exact, GW_MultiInit | 2317.2933 |
| 10×10 | Exact | 4525.7574 |
| 12×14 | GW_MultiInit | 5405.2937 |
| 15×12 | GW_MultiInit | 11089.4670 |
| 20×20 | GW_Default | 8405.1519 |
| 30×30 | GW_MultiInit | 13205.9421 |
| 40×50 | GW_MultiInit | 21181.4378 |
| 50×40 | GW_MultiInit | 27032.4820 |

| Instance (Agents×Tasks) | Best Method | Objective |
|---|---|---|
| 60×60 | GW_MultiInit, GW_Default | 28256.4427 |
| 100×100 | GW_MultiInit | 88468.94 |

Table 7 compares the objective values obtained by different Gromov-Wasserstein-based methods, a genetic algorithm (GA), and the exact solver across small and moderate CQAP instances where exact solutions are feasible. For each instance, the total mass is listed along with the objective values produced by:

- GW Default and GW MultiInit,
- Entropic GW (EGW) with three regularization parameters (ε = 0.3, 0.5, 0.7),
- Fused GW (FGW) with three blending weights (α = 0.3, 0.5, 0.7),
- Genetic Algorithm (GA), and
- The Exact method. The Exact method provides a baseline for comparison by solving the CQAP to global optimality using Gurobi's mixed-integer quadratic programming (MIQP) solver.
- The implementation uses Gurobi 9.x as the optimization engine, with the NonConvex=2 parameter enabling Gurobi's spatial branch-and-bound algorithm specifically designed for non-convex quadratic problems.

Bold values indicate the best result for each problem. GW MultiInit stays close to or matches the exact solution, making it the most reliable approximate method for all problem sizes. GA works well on small problems and often matches the exact result, but it gets less reliable as the problem grows. EGW and FGW are usually less accurate than GW MultiInit, but they give more flexibility and can handle large or noisy problems better. The key point is that GW MultiInit gives high-quality solutions while using much less computation than the exact method or GA. Exact methods stop working as problems get bigger (sometimes failing with interruptions). GA also takes much more time to compute. GW MultiInit stays efficient and keeps giving near-optimal results, so it is the best choice for practical use. Figure 3 shows the scalability of different methods. The exact solution provides the baseline, while other methods show varying degrees of approximation quality as problem size increases.

Table 7: Summary of Methods Performance (Problems with Exact Solutions)

| Problem Size | Mass | GW Default | GW MultiInit | EGW 0.3 | EGW 0.5 | EGW 0.7 | FGW 0.3 | FGW 0.5 | FGW 0.7 | GA | Exact |
|---|---|---|---|---|---|---|---|---|---|---|---|
| 3×3 | 9 | 995.89 | **981.63** | 995.89 | **981.6324** | 981.6325 | 1489.06 | **981.63** | **981.63** | **981.63** | **981.63** |
| 4×4 | 10 | **821.43** | **821.43** | **821.43** | 821.441 | 821.445 | 1180.88 | 1180.88 | 900.83 | **821.43** | **821.43** |

| Problem Size | Mass | GW Default | GW MultiInit | EGW 0.3 | EGW 0.5 | EGW 0.7 | FGW 0.3 | FGW 0.5 | FGW 0.7 | GA | Exact |
|---|---|---|---|---|---|---|---|---|---|---|---|
| 5×6 | 17 | 1430.52 | 1430.52 | 1430.53 | 1430.86 | 1432.81 | 2694.07 | 1622.19 | 1622.19 | 1432.81 | **1420.18** |
| 6×5 | 21 | 2324.62 | **2317.29** | 2324.62 | 2324.65 | 2325.10 | 3102.40 | 3102.40 | 3102.40 | 2324.62 | **2317.29** |
| 10×10 | 31 | 7729.06 | 4975.22 | 5653.20 | 5670.36 | 5752.01 | 5879.82 | 5364.28 | 5247.65 | 7729.06 | **4525.76** |
| 12×14 | 36 | 7319.04 | **5405.29** | 8339.81 | 7755.81 | 6347.12 | 8287.72 | 7149.91 | 5826.33 | 7319.04 | **5951.45*** |

*The best result until 4367.95 seconds after the run was interrupted

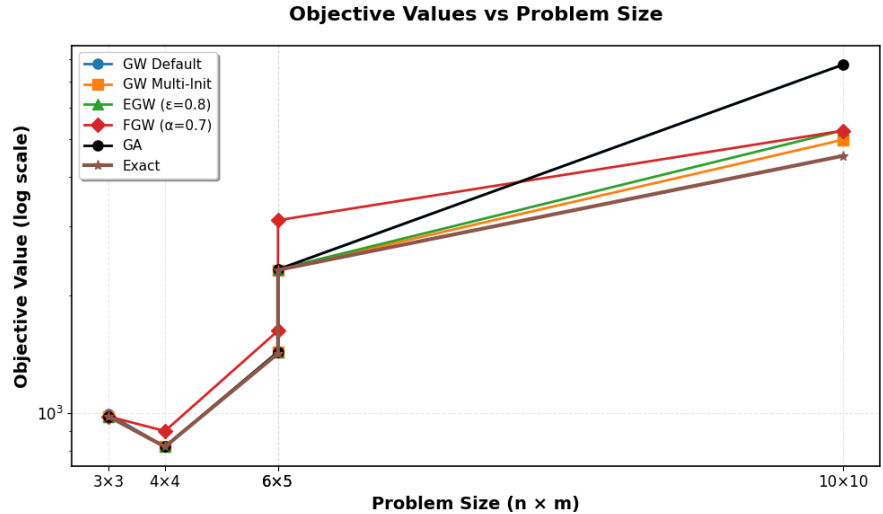

Figure 3: Runtime Scalability Analysis of GW Optimization Methods

Table 8 shows the time (in seconds) each method needs to solve CQAP instances of small to medium size, where exact solutions exist. The Exact method always gives the best solution, but its runtime grows fast as the problem gets bigger. For medium problems, it takes several minutes, and for larger ones, it becomes infeasible. The GA also shows significant computational overhead, requiring from a few seconds to over two minutes even for moderate-sized problems. In contrast, all Gromov-Wasserstein-based methods complete in under 2 seconds, even for the largest problems shown. These results show that GW methods scale better, especially when exact computation or GA takes too much time. Figure 4 shows how exact methods grow exponentially in runtime, while approximate methods keep the cost reasonable. Figure 5 reveals the trade-off between average solution quality and average computation time among test problems. Methods in the bottom-left represent the best compromise between speed and accuracy.

Table 8: Computation Times (seconds) for Exact-Solvable Problems

| Problem Size | GW Default | GW MultiInit | EGW 0.3 | EGW 0.5 | EGW 0.7 | FGW 0.3 | FGW 0.5 | FGW 0.7 | GA | Exact |
|---|---|---|---|---|---|---|---|---|---|---|
| 3×3 | **0.001** | 0.010 | 0.083 | 0.058 | 0.041 | **0.000** | 0.001 | **0.000** | 2.265 | 0.034 |
| 4×4 | **0.001** | 0.012 | 0.151 | 0.148 | 0.149 | **0.000** | 0.001 | **0.000** | 3.350 | 0.065 |

| Problem Size | GW Default | GW MultiInit | EGW 0.3 | EGW 0.5 | EGW 0.7 | FGW 0.3 | FGW 0.5 | FGW 0.7 | GA | Exact |
|---|---|---|---|---|---|---|---|---|---|---|
| 5×6 | **0.001** | 0.022 | 0.146 | 0.224 | 0.280 | **0.001** | **0.001** | **0.001** | 11.239 | 0.801 |
| 6×5 | **0.001** | 0.021 | 0.065 | 0.056 | 0.053 | **0.000** | **0.001** | **0.000** | 40.445 | 1.254 |
| 10×10 | 0.005 | 0.084 | 0.304 | 0.311 | 1.179 | **0.004** | **0.004** | **0.004** | 77.509 | 732.42 |
| 12×14 | 0.018 | 0.229 | 2.028 | 2.036 | 0.510 | 0.037 | 0.039 | **0.015** | 134.691 | 4367.95 |

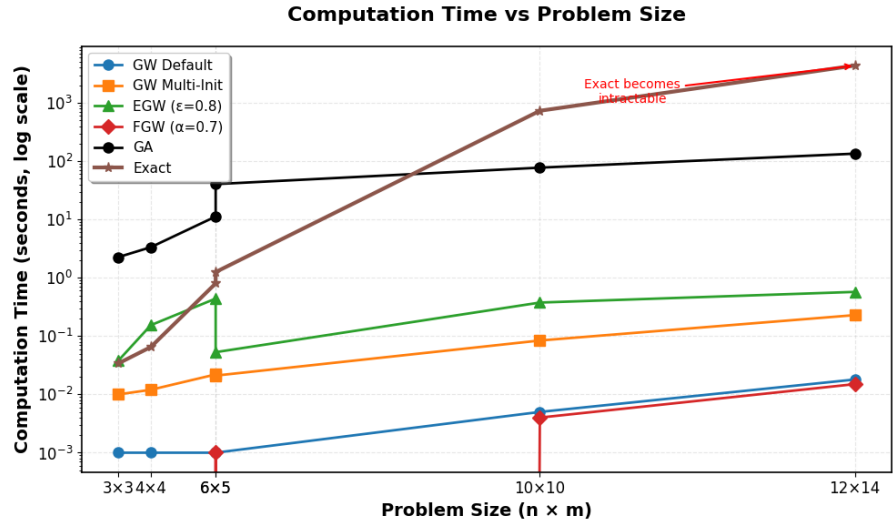

Figure 4: Computational Time Growth Comparison: Exact vs Approximate Methods

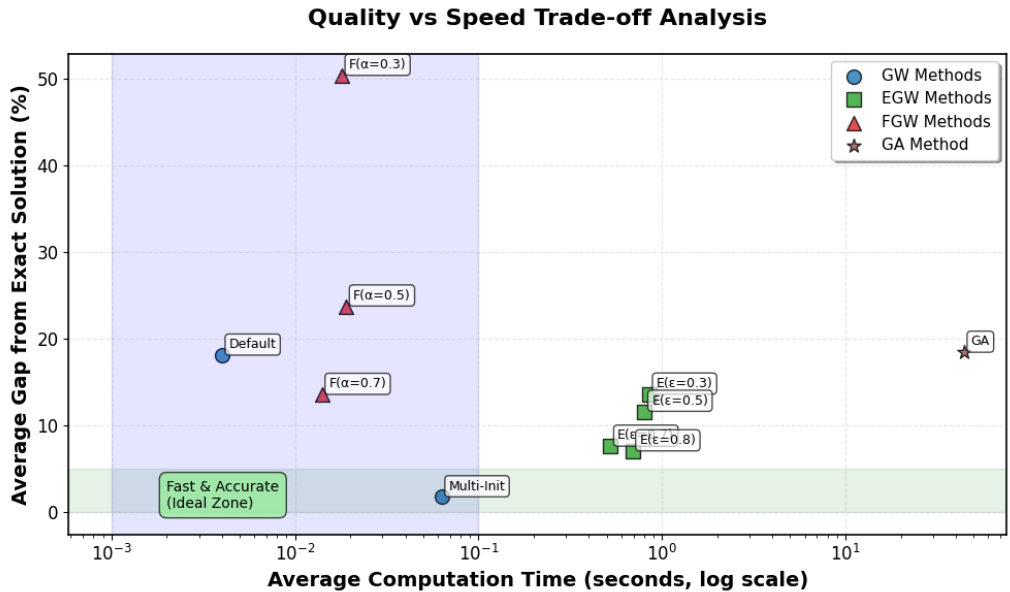

Figure 5: Trade-off Between Average Solution Quality and Average Computational Time

Table 9 reports the percentage gap between each method's objective value and the exact solution for problems where the exact result is known. The gap is computed as a relative error:

$$Gap\ (\%) = \frac{Approx.\ Objective - Exact\ Objective}{Exact\ Objective} \times 100$$

Lower values mean better accuracy. MultiInit_GW gives the smallest average gap (1.78%) and often matches the exact solution. EGW also performs well, with the gap getting smaller as the regularization strength (ε) grows. FGW shows larger gaps, especially with lower α values, which shows it is sensitive to how feature and structure costs are blended. Overall, the results show that GW MultiInit and EGW with higher ε give the best balance between accuracy and computational efficiency for approximate CQAP solutions.

Table 9: Gap from Exact Solution (%)

| Problem Size | GW Default | GW MultiInit | EGW 0.3 | EGW 0.5 | EGW 0.7 | FGW 0.3 | FGW 0.5 | FGW 0.7 | GA |
|---|---|---|---|---|---|---|---|---|---|
| 3×3 | 1.45 | 0.00 | 1.45 | 0.00 | 0.00 | 51.69 | 0.00 | 0.00 | 0.00 |
| 4×4 | 0.00 | 0.00 | 0.03 | 0.00 | 0.00 | 43.76 | 43.76 | 9.67 | 0.00 |
| 5×6 | 0.73 | 0.73 | 0.73 | 0.75 | 0.89 | 89.70 | 14.22 | 14.22 | 0.89 |
| 6×5 | 0.32 | 0.00 | 0.32 | 0.32 | 0.34 | 33.88 | 33.88 | 33.88 | 0.32 |
| 10×10 | 70.78 | 9.93 | 24.91 | 25.29 | 27.09 | 29.92 | 18.53 | 15.95 | 70.78 |
| 12×14 | 35.40 | 0.00 | 54.29 | 43.49 | 17.42 | 53.33 | 32.28 | 7.79 | 35.40 |
| **Average** | 18.12 | **1.78** | 13.62 | 11.64 | 7.62 | 50.38 | 23.78 | 13.59 | 0.32 |

Table 10 evaluates the effect of varying the entropic regularization parameter (ε) in the EGW method across several problem instances. The top half reports the gap from the exact solution (%) for each ε value. The bottom row shows the average runtime. Larger ε values (like 1.0 and 0.9) usually give smaller gaps, especially on bigger instances, but the runtime goes up a little. The best result is with ε = 0.8. It gives the lowest average gap (7.08%) and a moderate runtime (0.69 s). This means ε = 0.8 is a good choice. It balances accuracy and speed, so it can be used as the default for EGW in CQAP problems. Figure 6 shows that ε=0.8 provides the optimal balance for EGW methods, achieving the lowest average gap (7.08%) with reasonable computation time.

**Table 10: EGW Epsilon Parameter Analysis**

| Problem Size | ε = 1.0 | ε = 0.9 | ε = 0.8 | ε = 0.7 | ε = 0.6 | ε = 0.5 | ε = 0.3 |
|---|---|---|---|---|---|---|---|
| **Gap from Exact (%)** | | | | | | | |

| Problem Size | ε = 1.0 | ε = 0.9 | ε = 0.8 | ε = 0.7 | ε = 0.6 | ε = 0.5 | ε = 0.3 |
|---|---|---|---|---|---|---|---|
| 3×3 | - | - | - | 0.00 | 0.00 | 0.00 | 1.45 |
| 4×4 | 0.00 | 0.00 | 0.00 | 0.00 | 0.00 | 0.00 | 0.00 |
| 5×6 | 0.22 | 0.14 | 1.06 | 0.89 | 0.80 | 0.75 | 0.73 |
| 6×5 | 0.49 | 0.42 | 0.37 | 0.34 | 0.32 | 0.32 | 0.32 |
| 10×10 | 16.81 | 16.45 | 16.08 | 27.09 | 26.45 | 25.29 | 24.91 |
| 12×14 | 19.60 | 18.73 | 17.97 | 17.42 | 44.28 | 43.49 | 54.29 |
| **Average Gap** | 7.43 | 7.12 | **7.08** | 7.62 | 11.97 | 11.64 | 13.62 |
| **Computation Time (s)** | | | | | | | |
| Average | 0.63 | 0.64 | 0.69 | **0.52** | 0.79 | 0.80 | 0.85 |

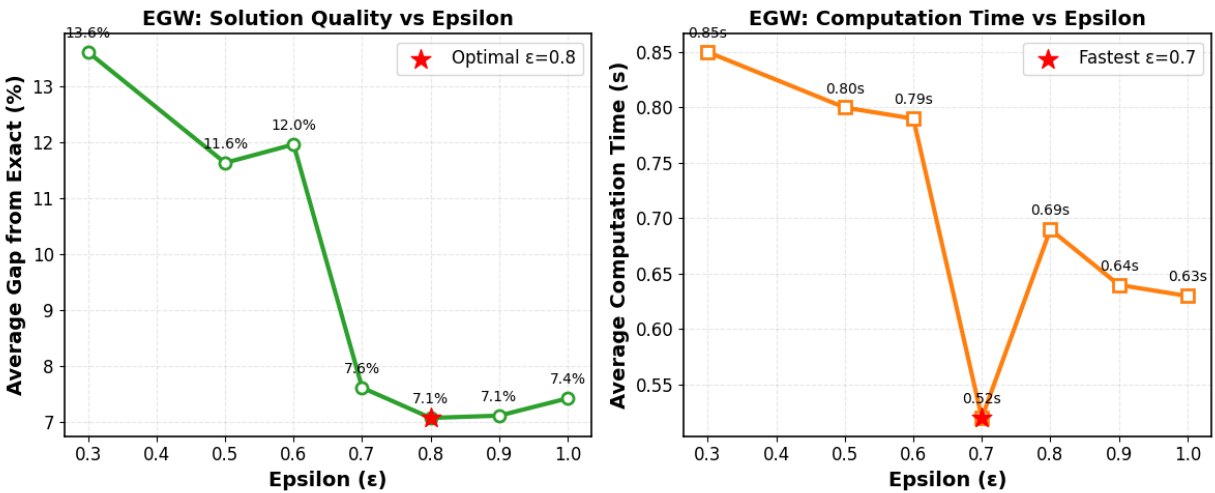

Figure 6: EGW Epsilon Parameter Optimization: Balancing Accuracy and Efficiency

Table 11 examines the impact of the α parameter in the FGW method, which controls the trade-off between structure preservation (Gromov-Wasserstein, α → 1) and feature alignment (classical Optimal Transport, α → 0). The top section shows the gap from the exact solution (%), while the bottom row reports the average computation time for each setting. The results clearly show that higher α values lead to significantly better accuracy. In particular, α = 0.7 achieves the lowest average gap (13.59%) with no added computational cost, suggesting that structure plays a critical role in solving CQAP effectively. Conversely, low α values (e.g., 0.0 and 0.3) result in large gaps, indicating that ignoring structural information greatly degrades solution quality. All configurations remain extremely fast (under 0.02 s), making FGW with α = 0.7 a strong choice when incorporating both structure and features.

Table 11: FGW Alpha Parameter Analysis

| Problem Size | α = 0.0 | α = 0.3 | α = 0.5 | α = 0.7 |
|---|---|---|---|---|
| **Gap from Exact (%)** | | | | |
| 3×3 | 51.69 | 51.69 | **0.00** | **0.00** |
| 4×4 | 43.76 | 43.76 | 43.76 | **9.67** |
| 5×6 | 89.42 | 89.70 | 14.22 | **14.22** |
| 6×5 | 32.14 | **33.88** | **33.88** | **33.88** |
| 10×10 | 33.97 | 29.92 | 18.53 | **15.95** |
| 12×14 | 66.80 | 53.33 | 32.28 | **7.79** |
| **Average Gap** | 52.96 | 50.38 | 23.78 | **13.59** |
| **Computation Time (s)** | | | | |
| Average | **0.014** | **0.018** | **0.019** | **0.014** |

Figure 7 demonstrates that α=0.7 provides the best performance for FGW methods, significantly improving solution quality while with the lowest average gap (13.59%), maintaining fast computation times.

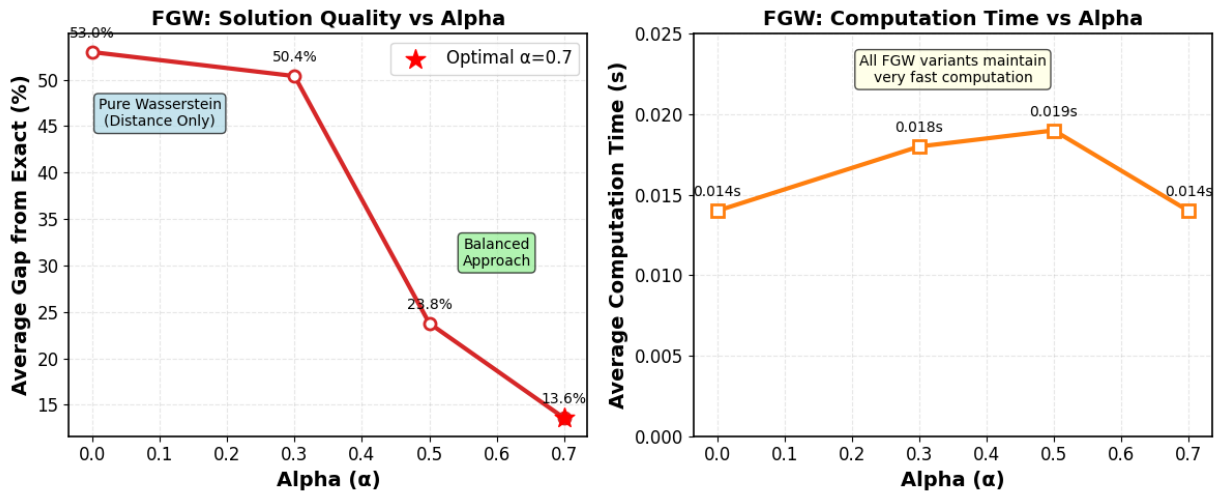

Figure 7: FGW Performance Analysis: Optimal Alpha Parameter Selection

Table 12 presents the objective values obtained by different methods on large-scale CQAP instances where exact solutions are unavailable due to computational limitations. Among all methods evaluated, GW MultiInit consistently achieves the best (lowest) objective values across all problem sizes, highlighting its robustness and scalability. The GA shows mixed performance, sometimes matching GW Default but generally falling short of GW MultiInit's superior results. While other methods such as EGW and FGW also produce feasible solutions, they often result in significantly higher objective values, especially in the largest instances. Notably, FGW's performance degrades sharply in high-dimensional settings, likely due to its sensitivity to feature alignment over structural consistency. These results confirm the suitability of GW MultiInit as a practical and effective solver for large CQAP problems where exact optimization is infeasible, outperforming both traditional metaheuristics like GA and other GW variants.

Table 12: Large Problems Performance (No Exact Solution Available)

| Problem Size | Mass | GW Default | GW MultiInit | EGW 0.3 | EGW 0.8 | FGW 0.3 | FGW 0.7 | GA | Best Method |
|---|---|---|---|---|---|---|---|---|---|
| 15×12 | 49 | 9657.404 | **8349.234** | 9411.960 | 9347.94 | 13574.982 | 14434.627 | 9657.404 | **GW MultiInit** |
| 20×20 | 62 | 8405.15 | **8405.15** | 12115.71 | 12737.47 | 8428.60 | 8471.25 | 8405.15 | **GW MultiInit** |
| 30×30 | 86 | 13304.89 | **13205.94** | 13572.06 | 15004.71 | 21060.28 | 19424.06 | 13304.89 | **GW MultiInit** |
| 40×50 | 118 | 22767.37 | **21181.44** | 23504.53 | 33807.78 | 34714.27 | 33806.34 | 22767.37 | **GW MultiInit** |
| 50×40 | 146 | 52157.45 | **27032.48** | 39587.11 | 49224.38 | 45991.30 | 45560.61 | 39587.11 | **GW MultiInit** |
| 60×60 | 176 | 28256.44 | **28256.44** | 30728.57 | 38284.84 | 41041.32 | 40649.58 | NA | **GW MultiInit** |
| 100×100 | 294 | 95867.89 | **88468.94** | 97709.38 | 128374.23 | 92070.45 | 89252.85 | NA | **GW MultiInit** |

Table 12 compares the computational efficiency of different methods on large CQAP instances. The FGW variant consistently demonstrates the fastest runtimes for medium-to-large problem sizes (20×20 and above), making it a practical choice when speed is critical. For the smallest large problem (15×12), both GW Default and FGW variants share the fastest runtime, finishing in just milliseconds (Table 13). Although GW MultiInit generally achieves better objective values, it requires substantially more computation time, often by an order of magnitude or more. EGW methods, while moderate in runtime, are still slower than FGW for large problems but demonstrate significantly better scaling characteristics compared to standard GW approaches. Looking at how the methods scale shows that EGW gets more competitive compared to GW as problem size grows. GW Default runtime goes from 0.013s to 57.83s (a 4,448× increase) from the smallest to the largest problem. EGW 0.3 increases more slowly, from 0.319s to 50.21s (only a 157× increase). GW MultiInit also scales poorly, with a 3,557× runtime increase, while EGW grows more modestly. EGW scales better because its Sinkhorn-based iterations handle dense matrix operations more efficiently than the conditional gradient approach, which repeatedly solves linear assignment subproblems. As problem size grows beyond what we tested, EGW's advantage over GW methods is expected to grow. The GA performs the worst. It takes much more time than any GW method, with runtimes from minutes to hours for large problems. For the largest tests (60×60 and 100×100), the GA didn't converge in a reasonable time, so it was left out of the comparison. This shows a clear trade-off between solution quality and speed. FGW 0.3 gives fast approximate solutions. GW MultiInit gives better quality at a moderate cost. GA can give competitive solutions but is very slow. Figure 8 shows that for large problems, GW MultiInit gives the best solution quality but takes much more computation time. EGW and FGW methods offer good compromises.

Table 13: Computation Time Analysis for Large Problems (seconds)

| Problem Size | GW Default | GW MultiInit | EGW 0.3 | EGW 0.8 | FGW 0.3 | FGW 0.7 | GA | Fastest Method |
|---|---|---|---|---|---|---|---|---|
| 15×12 | **0.013** | 0.238 | 0.319 | 2.55 | **0.013** | 0.013 | 293.332 | **GW Default/FGW** |
| 20×20 | 0.059 | 1.024 | 0.480 | 1.18 | **0.058** | 0.059 | 1283.793 | **FGW 0.3** |
| 30×30 | 0.294 | 5.114 | 1.173 | 1.39 | **0.289** | 0.292 | 4650.805 | **FGW 0.3** |
| 40×50 | 1.469 | 24.381 | 5.086 | 1.91 | **1.412** | 1.420 | 19168.123 | **FGW 0.3** |
| 50×40 | 1.442 | 24.349 | 2.164 | 2.07 | **1.418** | 1.429 | 25276.4331 | **FGW 0.3** |
| 60×60 | 4.653 | 78.373 | 5.358 | 5.29 | **4.590** | 4.603 | NA | **FGW 0.3** |
| 100×100 | 57.83 | 846.45 | 50.21 | 41.48 | 38.29 | **37.49** | NA | **FGW 0.7** |

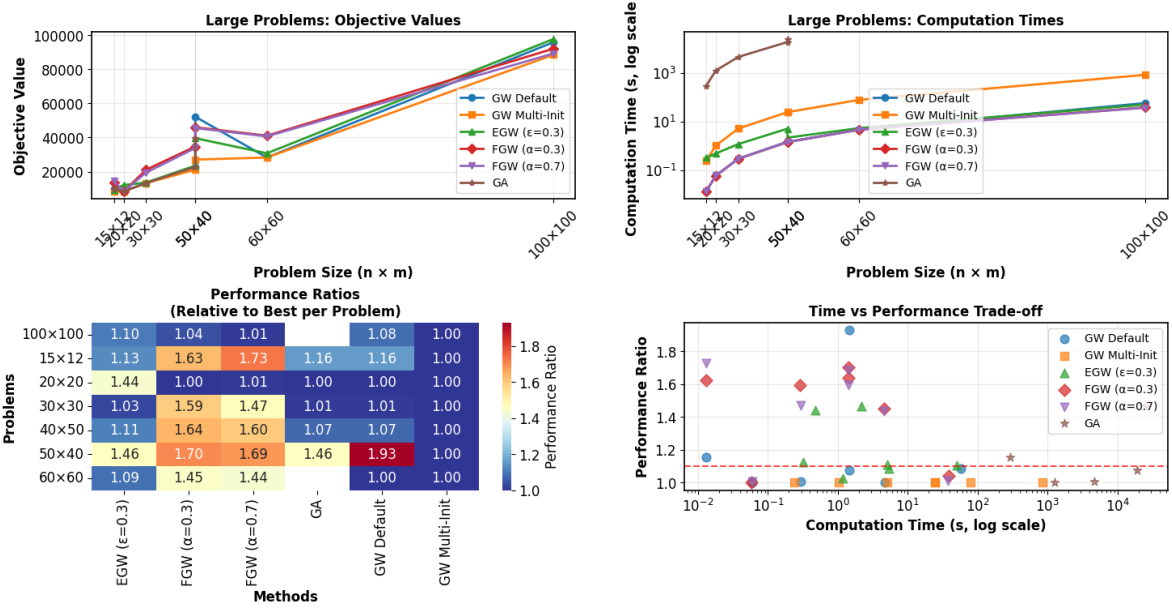

Figure 8: Solution Quality and Runtime Trade-offs for Large-Scale Problems

## 6. Conclusion, limitations, and perspectives

Gromov–Wasserstein (GW) optimal transport provides a powerful framework for comparing metric-measure spaces while preserving geometric structure, directly linking to the non-convex Capacitated Quadratic Assignment Problem (CQAP), which is NP-hard. In this paper, we investigated several GW variants and solvers, including entropy-regularized GW, Fused GW (FGW), and genetic algorithm (GA) approaches, and introduced multiple initialization GW (multi-int GW) as a more effective alternative. Gromov–Wasserstein (GW) optimal transport provides a powerful framework for comparing metric-measure spaces while preserving geometric structure, directly linking to the non-convex Capacitated Quadratic Assignment Problem (CQAP), which is NP-hard. In this paper, we investigated several GW variants and solvers, including entropy-regularized GW, Fused GW (FGW), and genetic algorithm (GA) approaches, and introduced multiple initialization GW (multi-int GW) as a more effective alternative. Our results show that multi-int GW consistently gets higher alignment accuracy and is more robust across synthetic benchmarks and real-world datasets. It performs better than EGW, FGW, and GA, and is still manageable to compute using Sinkhorn-based regularization. We looked at performance across benchmarks and found that approximate and regularized solvers can cut computational cost a lot while keeping high-quality solutions, especially for large CQAP instances.

Looking ahead, several directions could make GW methods more scalable and flexible. Some promising options are multi-marginal GW for aligning multiple metric-measure spaces at once, unbalanced GW for measures with unequal mass, partial GW for missing or noisy correspondences, sliced GW for computation that uses less dimension, and scalable GW for billion-scale graph and point cloud matching. Algorithm improvements like convex relaxations, spectral correspondence techniques, scalable graph-based GW learning, sampled GW, importance sparsification, and low-rank GW can also lower complexity without losing accuracy. Bayesian optimization in permutation spaces—with Gaussian process surrogates and kernels like Mallows or Merge—offers a solid way to deal with non-convexity in GW-related QAPs. This approach can find globally convergent strategies that work alongside existing gradient methods. These developments could let GW-based methods be used in larger and more complex areas, from structured data analysis and machine learning to network science and computational biology